\documentclass[a4paper,10pt]{article}

\newtheorem{thm}{Theorem}
\newtheorem{cor}[thm]{Corollary}

\newtheorem{lem}[thm]{Lemma}
\newtheorem{prop}[thm]{Proposition}

\def \b{{\beta}}

 \def \d{{\delta}}

\def \O{{\Omega}}
\def \p{{\varphi}}

\def \m{{\mu}}

\font\phh=cmcsc10  at  8 pt

\def \A{{\cal A}}

\def \P{{\bf P}}

\def \qq{{\qquad}}

\def \Z{{\bf Z}}

\def \e{{\varepsilon}}

\def \noi{{\noindent}}

 \def\qed{\hbox{\vrule height 6pt depth 0pt width
6pt}}
\def\cqfd{\hfill\penalty 500\kern 10pt\qed\medbreak}
 \font\phh=cmcsc10
at  8 pt
\scrollmode
\title{A Sharp Estimate for Divisors of Bernoulli Sums}
 \author{Michel Weber}

\begin{document}

\maketitle
%%%%%%%%%%%%%%%%%%%%%%%%%%%%%%%%%%%%%%%%%%%%%%%%%%%%%%%%%%%%%%%%%%%%%%%%%%%%%%%%%%%%%%%%%%%%%%

  \begin{abstract} \noi Let $S_n=\e_1+\ldots+\e_n$,  where   $ \e_i  $ are i.i.d.  Bernoulli r.v.'s.  Let $0\le r_d(n)<2d$ be   the
least residue of $n$  mod$(2d)$, $\bar r_d(n)= 2d -r_d(n)$ and $\b(n,d)=\max  ( {1\over d},  { 1\over \sqrt n}  )[e^{-   {r_d(n)^2/2 n}}
+e^{-   {\bar r_d(n)^2/2 n}}]$. We show that   
$$\sup_{2\le d\le n} \big|\P\big\{d|S_n\big\}- E(n,d)  \big|= {\cal O}\big({ \log^{5/2} n  \over n^{ 3/2}}\big),
$$  where $E(n,d) $ verifies $c_1\b(n,d)\le  E(n,d)\le c_2\b(n,d) $  and  $c_1,c_2 $ are numerical constants. 
   \end{abstract}

%%%%%%%%%%%%%%%%%%%%%%%%%%%%%%%%%%%%%%%%%%%%%%%%%%%%%%%%%%%%%%%%%%%%%%%%%%%%%%%%%%%%%%%%%%%%%%%%%%%%%%%%%%%%%%%%%%%
%%%%%%%%%%%%%%%%%%%%%%%%%%%%%%%%%%%%%%%%%%%%%%%%%%%%%%%%%%%%%%%%%%%%%%%%%%%%%%%%%%%%%%%%%%%%%%%%%%%%%%%%%%%%%%%%%%%
%%%%%%%%%%%%%%%%%%%%%%%%%%%%%%%%%%%%%%%%%%%%%%%%%%%%%%%%%%%%%%%%%%%%%%%%%%%%%%%%%%%%%%%%%%%%%%%%%%%%%%%%%%%%%%%%%%%
%%%%%%%%%%%%%%%%%%%%%%%%%%%%%%%%%%%%%%%%%%%%%%%%%%%%%%%%%%%%%%%%%%%%%%%%%%%%%%%%%%%%%%%%%%%%%%%%%%%%%%%%%%%%%%%%%%%
%%%%%%%%%%%%%%%%%%%%%%%%%%%%%%%%%%%%%%%%%%%%%%%%%%%%%%%%%%%%%%%%%%%%%%%%%%%%%%%%%%%%%%%%%%%%%%%%%%%%%%%%%%%%%%%%%%%
%%%%%%%%%%%%%%%%%%%%%%%%%%%%%%%%%%%%%%%%%%%%%%%%%%%%%%%%%%%%%%%%%%%%%%%%%%%%%%%%%%%%%%%%%%%%%%%%%%%%%%%%%%%%%%%%%%%
%%%%%%%%%%%%%%%%%%%%%%%%%%%%%%%%%%%%%%%%%%%%%%%%%%%%%%%%%%%%%%%%%%%%%%%%%%%%%%%%%%%%%%%%%%%%%%%%%%%%%%%%%%%%%%%%%%%
%%%%%%%%%%%%%%%%%%%%%%%%%%%%%%%%%%%%%%%%%%%%%%%%%%%%%%%%%%%%%%%%%%%%%%%%%%%%%%%%%%%%%%%%%%%%%%%%%%%%%%%%%%%%%%%%%%%
\section{Main result}
 Let    $\{\e_i, \,
i\ge 1\}$ denote a Bernoulli sequence defined on a joint probability space $(\tilde  \O,\tilde \A,\tilde \P)$,  with partial sums
$S_n=\e_1+\ldots+\e_n$.   Consider the    Theta function  
$$\displaystyle{ \Theta (d,m)  =  \sum_{\ell\in \Z} e^{im\pi{\ell\over   d }-{m\pi^2\ell^2\over 2 d^2}}. 
   } $$The improvment of the following result, which is    Theorem II in  \cite{We}, is the main purpose of this work.  
 \begin{prop}\label{Lemma4}      We have the following
uniform estimate: 
\begin{equation}\sup_{2\le d\le n}\Big|\P\big\{d|S_n\big\}- {\Theta(d,n)\over d}  \Big|= {\cal O}\big({ \log^{5/2} n  \over n^{ 3/2}}\big).
\end{equation} 
\end{prop}  
 This estimate is sharp already when $d<     ({ B  n/ \log n})^{1/2}$, otherwise 
   \begin{equation} 
  \big|\P\big\{d|S_n\big\}- {1\over d}  \big|\le  \cases{  C\big\{{ \log^{5/2} n  \over n^{ 3/2}}+ { 1\over d}
 e^{ - {n \pi^2 \over 2d^2}} \big\}   &  \qq if $d\le n^{ 1/2}$,\cr 
      {C \over n^{ 1/2}}   &  \qq if $ n^{ 1/2}\le d\le n$. 
}   \label{(2.10)}
\end{equation} 
   And this   is no longer efficient  when $d\gg \sqrt n$. The purpose of this Note is to remedy this by showing the existence of 
an extra corrective exponential factor in that case.
    Introduce a notation. Let
 $n\ge d\ge 2$ be integers and denote  by
 $r_d(n) $  the least residue of $n$ modulo $2d$:   
 $n\equiv r
 $ mod$(2d)$ and  
 $0\le r<2d$. Let also denote $\bar r_d(n)= 2d-r_d(n)$.

\begin{thm} \label{t1}    We have  
$$\sup_{2\le d\le n}\big|\P\big\{d|S_n\big\}- E(n,d)  \big|= {\cal O}\big({ \log^{5/2} n  \over n^{ 3/2}}\big).
$$ 
where $E(n,d)$ satisfies
 $$ {1\over 2\sqrt{2\pi}}  \le {E(n,d)\over  \max \big( {1\over 2d}, 
{ 1\over \sqrt n} \big)\big[e^{-   {r_d(n)^2\over 2 n}}+e^{-   {\bar r_d(n)^2\over 2 n}}\big]}\le  {32\over \sqrt{2\pi}} .
$$
    \end{thm}
 This   exponential   factor   is effective when
$\min( r_d(n),  \bar r_d(n) )\gg \sqrt n$. Its importance   is easily seen through   
 the following example. 
\smallskip \par 
Let $0<c<1$ and let $1\le \p_1(n)\le  c\p_2(n)$ be non-decreasing.
  Suppose $d$ is such that $2d\ge \sqrt n \p_2(n)$ with $r_d(n)$ large  so that $\sqrt{ 
n}\p_1(n)
\le r_d(n)
\le c\sqrt{  n}\p_2(n)$. Then $\bar r_d(n)\ge (1-c) \sqrt n \p_2(n)$ and so 
$$    E(n,d)\le {32\over \sqrt{2\pi n}}    
  \big[e^{-   {  \p_1^2(n)  \over 2  }}+e^{-   { (1-c)^2   \p_2^2(n) \over 2  }}\big]     .
$$
 Let $0< A_1\le A_2$. By taking $\p_i(n)= \sqrt{2A_i\log n}$,  $i=1,2$, we get 
$$    E(n,d)\le C\max \big(    
  n^{-1/2-A_1},n^{-1/2-(1-c)^2    A_2} \big)\ll n^{-1/2}     .
$$
Thus  we get a much better upper bound than    in  (\ref{(2.10)}).  The
proof uses  estimates for Theta functions, which are provided in the next Section.

%%%%%%%%%%%%%%%%%%%%%%%%%%%%%%%%%%%%%%%%%%%%%%%%%%%%%%%%%%%%%%%%%%%%%%%%%%%%%%%%%%%%%%%%%%%%%%%%%%%%%%%%%%%%%%%%%%%
%%%%%%%%%%%%%%%%%%%%%%%%%%%%%%%%%%%%%%%%%%%%%%%%%%%%%%%%%%%%%%%%%%%%%%%%%%%%%%%%%%%%%%%%%%%%%%%%%%%%%%%%%%%%%%%%%%%
%%%%%%%%%%%%%%%%%%%%%%%%%%%%%%%%%%%%%%%%%%%%%%%%%%%%%%%%%%%%%%%%%%%%%%%%%%%%%%%%%%%%%%%%%%%%%%%%%%%%%%%%%%%%%%%%%%%
%%%%%%%%%%%%%%%%%%%%%%%%%%%%%%%%%%%%%%%%%%%%%%%%%%%%%%%%%%%%%%%%%%%%%%%%%%%%%%%%%%%%%%%%%%%%%%%%%%%%%%%%%%%%%%%%%%%
%%%%%%%%%%%%%%%%%%%%%%%%%%%%%%%%%%%%%%%%%%%%%%%%%%%%%%%%%%%%%%%%%%%%%%%%%%%%%%%%%%%%%%%%%%%%%%%%%%%%%%%%%%%%%%%%%%%
%%%%%%%%%%%%%%%%%%%%%%%%%%%%%%%%%%%%%%%%%%%%%%%%%%%%%%%%%%%%%%%%%%%%%%%%%%%%%%%%%%%%%%%%%%%%%%%%%%%%%%%%%%%%%%%%%%%
%%%%%%%%%%%%%%%%%%%%%%%%%%%%%%%%%%%%%%%%%%%%%%%%%%%%%%%%%%%%%%%%%%%%%%%%%%%%%%%%%%%%%%%%%%%%%%%%%%%%%%%%%%%%%%%%%%%

\section{Theta Function Estimates}

Let $E(n,d): = {\Theta(d,n)\over d} $.  By the Poisson summation formula  
$$\sum_{\ell\in \Z} e^{-(\ell+\d)^2\pi x^{-1}}=x^{1/2} \sum_{\ell\in \Z}  e^{2i\pi \ell\d -\ell^2\pi x}, 
$$ 
  where $x$  is any real and $0\le \d\le 1$, we get
  with the choices $x=\pi
n/(2d^2)$,
$\d=n/(2d)$ 
\begin{equation}\label{E}E(n,d)  =\sqrt{  { 2\over \pi
n}}\sum_{h\in \Z} e^{-2(\{ {n\over 2d}\} +h)^2 {d^2\over n}}.
\end{equation}
Let $a>0$, $0\le \m\le 1$ and write $\bar \m:=1-\m$.    We begin with elementary estimates of  
$$S(\m, a):= \sum_{h\in \Z} e^{-a ( \m +h)^2} = e^{-   a\m ^2}+ e^{- a\bar \m^2 }+  \sum_{h=1}^\infty e^{- a( \m +h)^2}+ \sum_{h=1}^\infty e^{- a( h+\bar \m
)^2} .$$

  \begin{lem} Define for    $0\le \m\le 1$ and $a>0$,
 $  \p(\m, a)={1\over \sqrt {2a} +2a\m
 }$.  Then   
$$( \p(\m, a)-1) e^{- a \m^2  } \le \sum_{h=1}^\infty e^{- a( \m +h)^2}\le 2\p(\m, a) e^{- a \m^2  }.$$
 \end{lem}
\noi{\it Proof.}   Consider  Mill's ratio  
 $R(x) = e^{x^2/2}\int_x^\infty e^{-t^2/2}\ dt$. Then (\cite{Mi} section 2.26)
  $${1\over 1 + x}\le  {2\over \sqrt{x^2+4} +x}\le R(x) \le 
{2\over \sqrt{x^2+{8/ \pi}} +x}\le 
{2\over 1+x},\qq x\ge 0.  $$ 
     First
$$ \int_{0 }^\infty e^{- a( \m +x)^2} dx-   e^{- a  \m  ^2}  \le
\sum_{h=1}^\infty e^{- a( \m +h)^2}\le  \int_{0 }^\infty e^{- a( \m +x)^2} dx.$$
 But 
 $$\int_{0 }^\infty e^{- a( \m +x)^2} dx = {e^{-  \m^2 a }\over \sqrt {2a}} R(\m \sqrt {2a}) \quad {\rm and}\quad {1\over 1 +  \m \sqrt {2a}} \le    
  R(\m \sqrt {2a}) \le  {2\over 1 +\m
\sqrt {2a}}. $$
  Thus 
$$  {1\over \sqrt {2a} +2a\m
 } e^{- a \m^2  }\le \int_{0 }^\infty e^{- a( \m +x)^2} dx\le  {2\over \sqrt {2a} +2a\m
 } e^{- a \m^2  } . $$
Hence 
$$  (\p(\m, a)-1)   e^{- a  \m  ^2} \le
\sum_{h=1}^\infty e^{- a( \m +h)^2}\le  2 \p(\m, a) e^{- a \m^2  },$$
as claimed. \cqfd 
    \begin{cor} \label{c1}Put 
 $  \psi(\m, a):=\big( 1+ 
  \p(\m, a) \big)e^{-   a\m ^2}$. Then for every $0\le \m\le 1$ and $  a>0$ 
$${1\over 2} \le { S(\m, a) \over \psi(\m, a)+ \psi(\bar \m, a)}\le 2   .$$
 
\end{cor}
{\it Proof.}   At first by the previous Lemma
$$ A :=  e^{-   a\m ^2} +  \sum_{h=1}^\infty e^{- a( \m +h)^2} \le  e^{-   a\m ^2} \big(1+ 2\p(\m, a)\big).
$$
   Next 
$$  A\ge e^{-   a\m ^2} +  {1\over 2}\sum_{h=1}^\infty e^{- a( \m +h)^2} \ge 
e^{-   a\m ^2} + 
{1\over 2} \big( \p(\m, a)- 1  \big)  e^{- a  \m  ^2}={1\over 2}\big( 1+ 
  \p(\m, a) \big)e^{-   a\m
^2}  .$$
Thereby
 ${1/2}\le {A\over\psi(\m, a)}\le 2$. 
 Operating similarly with $\bar A =e^{-   a\bar\m ^2} +  \sum_{h=1}^\infty e^{- a( \bar\m +h)^2}$   leads to
$${1\over 2}\le {S(\m, a)\over \psi(\m, a)+ \psi(\bar \m, a)}\le 2.$$
\cqfd

 Notice  that $\p(0, a)={1/\sqrt {2a}}$  and   
\begin{equation}\label{0}{1\over 2}\big( 1+ 
{1\over \sqrt {2a}} \big)\ \le S(0, a)= 1+2\sum_{h=1}^\infty e^{-    a h ^2  } \le 4\big( 1+ 
{1\over \sqrt {2a}} \big).
\end{equation}

We now need an extra Lemma. 
\begin{lem} \label{l3} Let $n=2d K    +r$ with $1\le r\le 2d$. Then 
$$   {1\over 2} \max\big( {1\over 2d}, {1\over \sqrt n}\big)  e^{-   {r^2\over 2 n}}   \le  {\psi({r\over 2d}, {2d^2\over n})  \over \sqrt
n} 
\le  
   2\max\big( {1\over 2d}, {1\over \sqrt n}\big)  e^{-   {r^2\over 2 n}} . $$

\end{lem}
{\it Proof.}   We have 
$$  \psi({r\over 2d}, {2d^2\over n})    =\Big( 1+ 
  {\sqrt n\over  2d } \,{1\over 1 + {r\over \sqrt n}}  \Big)e^{-   {r^2\over 2n}}.  $$
 We
consider three cases. \smallskip\par \noi {\bf Case a.} $2d\le \sqrt n$. Then ${r\over \sqrt n}< {2d\over \sqrt n}\le 1$, and so 
  $  
  {\sqrt n\over 4d }    e^{-   {r^2\over 2 n}}\le  \psi({r\over 2d}, {2d^2\over n})  \le  
  {\sqrt n\over  d }  e^{-   {r^2\over 2 n}} $, 
which implies 
 $$ 
{1\over 2}  \max\big( {1\over 2d}, {1\over \sqrt n}\big)  e^{-   {r^2\over 2 n}}=     {e^{-   {r^2\over 2 n}}\over  4  d  }     \le 
{\psi({r\over 2d}, {2d^2\over n})  \over
\sqrt n} 
\le  
  {e^{-   {r^2\over 2 n}}\over  d }    =2\max\big( {1\over 2d}, {1\over \sqrt n}\big)  e^{-   {r^2\over 2 n}}. $$
  {\bf Case b.}   $2d\ge \sqrt n$ and $r \le \sqrt n $. Here  we have
  $ e^{-   {r^2\over 2 n}}\le  \psi({r\over 2d}, {2d^2\over n})  \le 2e^{-   {r^2\over 2 n}} $,
which implies 
  $$    \max\big( {1\over 2d}, {1\over \sqrt n}\big)  e^{-   {r^2\over 2 n}}=
     {e^{-   {r^2\over 2 n}}\over   \sqrt{  n } }     \le  {\psi({r\over 2d}, {2d^2\over n})  \over \sqrt n} 
\le  
  {2 e^{-   {r^2\over 2 n}}\over  \sqrt n }  = 2\max\big( {1\over 2d}, {1\over \sqrt n}\big)  e^{-   {r^2\over 2 n}} . $$
  {\bf Case c.}   $2d\ge \sqrt n$ and $r\ge  \sqrt n $. The  exponential factor $e^{-   {r^2\over 2 n}}$  
may this time play  a role (if $r\gg \sqrt n$), and we have  
  $ e^{-   {r^2\over 2 n}}\le  \psi({r\over 2d}, {2d^2\over n})  \le
{3\over 2}e^{-   {r^2\over 2 n}}$  which implies 
  $$    \max\big( {1\over 2d}, {1\over \sqrt n}\big)  e^{-   {r^2\over 2 n}}=
     { e^{-   {r^2\over 2 n}}\over   \sqrt{  n } }    \le  {\psi({r\over 2d}, {2d^2\over n})  \over \sqrt n} 
\le  
  {3 e^{-   {r^2\over 2 n}}\over  2\sqrt n }  = {3\over 2}\max\big( {1\over 2d}, {1\over \sqrt n}\big)  e^{-   {r^2\over 2 n}} . $$
 
Summarizing cases a) to c), we  have that
$$   {1\over 2} \max\big( {1\over 2d}, {1\over \sqrt n}\big)  e^{-   {r^2\over 2 n}}   \le  {\psi({r\over 2d}, {2d^2\over n})  \over \sqrt
n} 
\le  
   2\max\big( {1\over 2d}, {1\over \sqrt n}\big)  e^{-   {r^2\over 2 n}} . $$
  \cqfd 
%%%%%%%%%%%%%%%%%%%%%%%%%%%%%%%%%%%%%%%%%%%%%%%%%%%%%%%%%%%%%%%%%%%%%%%%%%%%%%%%%%%%%%%%%%%%%%%%%%%%%%%%%%%%%%%%%%%
%%%%%%%%%%%%%%%%%%%%%%%%%%%%%%%%%%%%%%%%%%%%%%%%%%%%%%%%%%%%%%%%%%%%%%%%%%%%%%%%%%%%%%%%%%%%%%%%%%%%%%%%%%%%%%%%%%%
%%%%%%%%%%%%%%%%%%%%%%%%%%%%%%%%%%%%%%%%%%%%%%%%%%%%%%%%%%%%%%%%%%%%%%%%%%%%%%%%%%%%%%%%%%%%%%%%%%%%%%%%%%%%%%%%%%%
%%%%%%%%%%%%%%%%%%%%%%%%%%%%%%%%%%%%%%%%%%%%%%%%%%%%%%%%%%%%%%%%%%%%%%%%%%%%%%%%%%%%%%%%%%%%%%%%%%%%%%%%%%%%%%%%%%%
%%%%%%%%%%%%%%%%%%%%%%%%%%%%%%%%%%%%%%%%%%%%%%%%%%%%%%%%%%%%%%%%%%%%%%%%%%%%%%%%%%%%%%%%%%%%%%%%%%%%%%%%%%%%%%%%%%%
%%%%%%%%%%%%%%%%%%%%%%%%%%%%%%%%%%%%%%%%%%%%%%%%%%%%%%%%%%%%%%%%%%%%%%%%%%%%%%%%%%%%%%%%%%%%%%%%%%%%%%%%%%%%%%%%%%%
%%%%%%%%%%%%%%%%%%%%%%%%%%%%%%%%%%%%%%%%%%%%%%%%%%%%%%%%%%%%%%%%%%%%%%%%%%%%%%%%%%%%%%%%%%%%%%%%%%%%%%%%%%%%%%%%%%%
 \section{Proof}
    A first case is simple. \smallskip\par \noi {\bf Case I.} $2d|n$.   We have 
 $E(n,d)   =\sqrt{  { 2\over \pi
n}}S(0 , {2d^2\over n})$.   But by (\ref{0})
 $${1\over 2}\max\big( 1, 
{\sqrt n \over 2d}\big)\le {1\over 2}\big( 1+ 
{\sqrt n \over 2d} \big) \le S(0, {2d^2\over n})   \le 4\big( 1+ 
{\sqrt n \over 2d} \big)\le 8\max\big( 1, 
{\sqrt n \over 2d}\big).
$$
 Hence
 \begin{equation}\label{rho=0}{1\over \sqrt{2\pi}} \max\big( {1\over 2d}, {1\over \sqrt n}\big)  \le E(n,d) \le {16\over \sqrt{2\pi}}
\max\big( {1\over 2d}, {1\over
\sqrt n}\big). 
\end{equation}
 \bigskip\par    
\smallskip\par \noi {\bf Case II.} Now if $2d\not| n$,  write $n= 2dK  +\rho$  with $0< \rho<2d$. In our setting $a={2d^2\over n}$,
$\m=\{ {n\over 2d}\}={\rho\over 2d}$ and by (\ref{E}), $E(n,d)=\sqrt{  { 2/ \pi
n}}\,S(\{ {n\over 2d}\}, {2d^2\over n})$.    Applying Lemma \ref{l3} with $r=\rho$ gives 
\begin{equation}   {1\over 2} \max\big( {1\over 2d}, {1\over \sqrt n}\big)  e^{-   { \rho^2\over 2 n}}   \le  {\psi({  \rho\over
2d}, {2d^2\over n}) 
\over
\sqrt n} 
\le  
   2\max\big( {1\over 2d}, {1\over \sqrt n}\big)  e^{-   {  \rho^2\over 2 n}} . 
\end{equation}

As to   $\psi(\bar \m, {2d^2\over n})  $, we have
  $\bar \m={2d-\rho\over 2d}:= {\bar \rho\over 2d}$ and   $0< \bar \rho<2d$. Applying Lemma \ref{l3} with $r=\bar \rho$ gives
\begin{equation}   {1\over 2} \max\big( {1\over 2d}, {1\over \sqrt n}\big)  e^{-   {\bar \rho^2\over 2 n}}   \le  {\psi({\bar \rho\over
2d}, {2d^2\over n}) 
\over
\sqrt n} 
\le  
   2\max\big( {1\over 2d}, {1\over \sqrt n}\big)  e^{-   {\bar \rho^2\over 2 n}} . 
\end{equation}
  Consequently, by Corollary \ref{c1}
   \begin{equation} \label{rho}  {1\over 2\sqrt{2\pi}}   \le {E(n,d)\over \max \big( {1\over 2d}, 
{ 1\over \sqrt n} \big)\big[e^{-   {\rho^2\over 2 n}}+e^{-   {\bar \rho^2\over 2 n}}\big]}\le  {8\over \sqrt{2\pi}} .
\end{equation}
When $\rho=0$, it follows from estimate (\ref{rho=0}) that
 \begin{eqnarray*}  { \max \big( {1\over
2d},  { 1\over \sqrt n} \big)\over  \sqrt{2\pi}} \big[{1+e^{-   {\bar \rho^2\over 2 n}}\over 2}\big]&\le& { \max \big( {1\over
2d},  { 1\over \sqrt n} \big)\over \sqrt{2\pi}}
    \le E(n,d)  \le  {16 \max \big( {1\over
2d},  { 1\over \sqrt n} \big)\over \sqrt{2\pi}}   
\cr &\le & {32 \max \big( {1\over
2d},  { 1\over \sqrt n} \big)\over \sqrt{2\pi}}  \big[{1+e^{-   {\bar \rho^2\over 2 n}}\over 2}\big]. 
\end{eqnarray*}
Finally in either case
 \begin{equation} \label{rho1}  {1\over 2\sqrt{2\pi}}  \le {E(n,d)\over  \max \big( {1\over 2d}, 
{ 1\over \sqrt n} \big)\big[e^{-   {\rho^2\over 2 n}}+e^{-   {\bar \rho^2\over 2 n}}\big]}\le  {32\over \sqrt{2\pi}} .
\end{equation}
 \cqfd
 
{ 
 \medskip\par\noi  
\noi {\phh Michel Weber:  IRMA, Universit\'e
Louis-Pasteur et C.N.R.S.,   7  rue Ren\'e Descartes, 67084
Strasbourg Cedex, France. \noi E-mail: \  \tt
weber@math.u-strasbg.fr}

\end{document}